\documentclass[12pt,fleqn]{article}

\usepackage{mathrsfs}
\usepackage{graphicx}
\usepackage{amsmath}
\usepackage{amsfonts}
\usepackage{amssymb}

\newtheorem{theorem}{Theorem}[section]
\newtheorem{lemma}[theorem]{Lemma}

\newtheorem{corollary}[theorem]{Corollary}
\numberwithin{equation}{section}

\begin{document}

\title{ Spectral
characterization of a specific class of trees \thanks{This research
was partially supported by the NSF of China(No.10571077).}}

\author{ Xiaoxia Fan, Yanfeng Luo\\
{\small Department of Mathematics, Lanzhou University,}\\
{\small Lanzhou, Gansu 730000, PR China}\\
{\small fanxx06@lzu.cn}}
\date{}

\maketitle

\begin{abstract}
In this paper, it is shown that the graph $T_4(p,q,r)$ is determined
by its Laplacian spectrum and there are no two non-isomorphic such
graphs which are cospectral with respect to adjacency spectrum.

\vskip 0.05in


\end{abstract}

\section{ Introduction }
Graphs considered in this paper are undirected graphs without loops
and multiple edges. Let $G$ be a simple graph with $n$ vertices.
Denote by $A(G)$ and $D(G)$ the adjacency matrix and the diagonal
matrix with the vertex degrees of $G$ on the diagonal, respectively.
The matrix $L(G) = D(G) -A(G)$ is called the {\it Laplacian matrix}
of $G$. Denote by $P(G,\lambda)$ the adjacency polynomial
$det(\lambda I-A(G))$ of $G$. The multiset of eigenvalues of $A(G)$
(resp., $L(G)$) is called the {\it adjacency} (resp., {\it
Laplacian}) {\it spectrum} of $G$. Since $A(G)$ and $L(G)$ are real
symmetric matrices, their eigenvalues are real numbers. So we can
assume that $\lambda_1 \geq \lambda_{2}\geq \cdots \geq \lambda_{n}$
and $\mu_{1} \geq \mu_{2} \geq \cdots \geq \mu_{n}$ are the
adjacency eigenvalues and the Laplacian eigenvalues of $G$,
respectively. Two graphs are said to be {\it cospectral} with
respect to the adjacency (resp. Laplacian) spectrum if they have the
same adjacency (resp. Laplacian) spectrum. A graph is said to be
{\it determined by its adjacency (resp., Laplacian) spectrum} if
there is no other non-isomorphic graph with the same adjacency
(resp., Laplacian) spectrum.

Determining what kinds of graphs are determined is an old problem,
which is far from resolved, in the theory of graph spectra. In their
paper \cite{van Dam}, the authors conjectured that almost all graphs
are determined by their spectrum. However, it seems hard to prove a
graph to be determined by its spectrum and only a few graphs have
been proved to be determined by their spectrum. Therefore it would
be interesting to find more examples of graphs which are determined
by their spectrum. For the background on this problem and related
topics, the reader can consult \cite{van Dam, Edwin}. For more
recent results which have not been cited in \cite{van Dam, Edwin},
we refer to \cite{Ro, Oimidi, GR, L.Shen, F} and their references
for details.

Because the problem above is very hard to deal with, van Dam and
Haemers \cite{van Dam} suggested a modest problem, say, ``which
trees are determined by their spectrum?'' This paper will give a
complete answer to this modified problem for a class of special
trees.

As usual, we denote by $P_k$ the path with $k$ vertices. Let $G$ be
a graph. Denote by ${\mathcal L}(G)$ the line graph of $G$. We
denote by $T_4(p,q,r)$ the graph shown in Fig.~1. $T_4(p,q,r)$ is a
tree with $4$ vertices of degree $3$. For a $T_4(p,q,r)$ graph, we
always assume that $1\leq p\leq q \leq r$. The reader is referred to
\cite{Bondy} for any undefined notion and terminology on graphs in
this paper.

\begin{figure}[htbp]
\begin{center}
\includegraphics[totalheight=5cm]{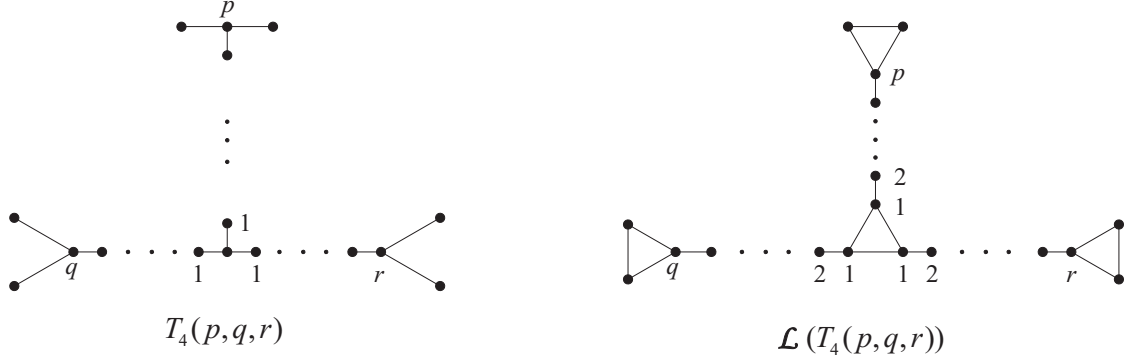}\\
 \caption{\label{T(pqr)} \small{The  graphs  $T_4(p,q,r)$ and  ${\mathcal L}(T_4(p,q,r))$ where $p,q,r\geq 1$. }}
\end{center}
\end{figure}

In this paper we will show that $T_4(p,q,r)$ is determined by its
Laplacian spectrum and there are no two non-isomorphic graphs which
are cospectral with respect to adjacency spectrum.

\section{ Preliminaries }

In this section, we will present some known results which will be used in this paper.

\begin{lemma}(\cite{Bondy}) \label{21}
Two trees $T$ and $T^{'}$  are cospectral with
respect to the Laplacian matrix if and only if their line graphs are
cospectral with respect to the adjacency matrix.
\end{lemma}

\begin{lemma}(\cite{Godsil})\label{20}
If ${\mathcal L}(G)\cong {\mathcal L}(H)$ with $\{G,H\}\neq
\{K_3,K_{1,3}\}$. Then $G\cong H$.
\end{lemma}

Let $W_n$ be the graph obtained from the path $P_ {n-2}$ (indexed in
natural order $1,2,\dots,n-2$) by adding two pendant edges at
vertices $2$ and $n-3$.

\begin{lemma}(\cite{Hoffman}) \label{uv}
Let $G$ be a connected graph that is not isomorphic to $W_{n}$ and $G_{uv}$ be the graph obtained from $G$ by
subdividing the edge $uv$ of $G$. If $uv$ lies on an internal path
of G, then $\lambda_{1}(G_{uv})\leq \lambda_{1}(G)$.
\end{lemma}

\begin{lemma}(\cite{van Dam}) \label{211}
Let $G$ be a graph. The following can be obtained from the adjacency spectrum and from the Laplacian
spectrum:

(i) The number of vertices,

(ii) The number of edges.

The spectrum of the adjacency matrix determines:

(iii) The number of closed walks of any length.

The Laplacian spectrum determines:

(iv) The number of spanning trees,

(v) The number of components,

(vi) The sum of squares of degrees of vertices.
\end{lemma}

Let $N_{G}(H)$ be the number of subgraphs of a graph G which is
isomorphic to $H$ and let $N_{G}(i)$ be the number of closed walks
of length $i$ of $G$.

\vskip 0.1in

\begin{lemma}(\cite{Oimidi}) \label{212} Let $G$ be a graph. Then

(i) $N_{G}(2) = 2m, \ N_{G}(3) = 6N_{G}(K_{3})$;

(ii) $N_{G}(4) = 2m + 4N_{G}(P_{3}) +
8N_{G}(C_{4})$, $N_{G}(5) = 30N_{G}(K_{3}) + 10N_{G}(C_{5}) +
10N_{G}(G_{1})$;

(iii) $ N_{G}(7) =
126N_{G}(K_{3})+84N_{G}(G_{1})+14N_{G}(G_{2})+14N_{G}(G_{3})+14N_{G}(G_{4})+28N_{G}(G_{5})+
42N_{G}(G_{6}) + 28N_{G}(G_{7}) + 112N_{G}(G_{8}) + 70N_{G}(C_{5}) +
14N_{G}(C_{7})$. (see Fig.~2).

\end{lemma}

\begin{figure}[htbp]
\begin{center}
\includegraphics[totalheight=6cm]{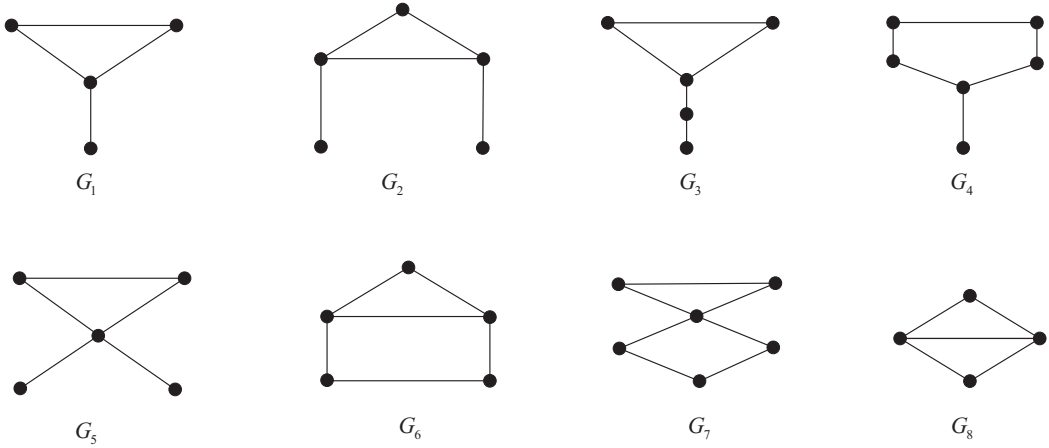}\\
 \caption{\label{lem} \small{The graphs $G_i$, i=1,\dots,8. }}
\end{center}
\end{figure}

\begin{lemma}(\cite{Li}) \label{213} Let $G$ be a graph with $V (G) \neq \emptyset $ and $E(G)\neq \emptyset$.
 Then $$\Delta(G) + 1 \leq \mu_{1}
\leq max \{ \frac {d_{u}(d_{u}+m_{u})+d_{v}(d_{v}+m_{v} )}
{d_{u}+d_{v} },uv \in E(G) \}$$ where $\Delta(G)$ denote the maximum
vertex degree of $G$, and $m_{v}$ the average of degrees of the
vertices adjacent to the vertex $v$ in $G$.
\end{lemma}

\begin{lemma}(\cite{Cvetkovic}, \cite{AJ})
\label{24}  Let $v$ be a vertex of a graph $G$ and let $C(v)$ denote
the collection of cycles containing $v$. Then the characteristic
polynomial of $G$ satisfies
 \begin{equation*} P(G,\lambda)=\lambda
P(G\setminus\{v\},\lambda) -\sum \limits_{u\thicksim v}
P(G\setminus\{u,v\},\lambda)-2\sum \limits_{Z\in C(v)} P(G\setminus
V(Z),\lambda).
\end{equation*}
\end{lemma}

For the sake of convenience, denote $P(P_r,\lambda)$ by
$p_r=p_r(\lambda)$. For convenience's sake, let $p_0=1, p_{-1}=0$
and $p_{-2}=-1$.

\begin{lemma}(\cite{F})\label{25} $p_r=\frac {x^{2r+2}-1} {x^{r+2}-x^{r}}$ and $p_r(2)=r+1$, where $x$ satisfies $x^{2}-\lambda x+1=0$.
\end{lemma}

A centipede is a graph obtained by appending a pendant vertex to
each vertex of degree $2$ of a path.

\begin{lemma}(\cite{Ro})\label{26}
The centipede is determined by its Laplacian spectrum.
\end{lemma}

\begin{lemma}(\cite{Edwin})\label{27} For bipartite graphs, the sum of cubes of degrees is determined by the Laplacian
spectrum.
\end{lemma}

\section{$T_4(p,q,r)$ is determined by its Laplacian spectra }

In this section, we will show that $T_4(p,q,r)$
is determined by its Laplacian spectrum. To this aim, we need to compute the
characteristic polynomial of the line graph ${\mathcal L}(T_4(p,q,r)$ of $T_4(p,q,r)$. By using Lemma \ref{24}
with $v$ being the vertices of degree three, we have
\begin{eqnarray*}
P({\mathcal L}(T_4(p,q,r)),\lambda) & =& f(q,r)( \lambda
h_{p-1}-h_{p-2})-h_{p-1}(h_{q-1}h_{r}-h_qh_{r-1}-2h_{q-1}h_{r-1}),\\
P({\mathcal L}(T_4(1,q,r)),\lambda) & =& f(q,r)( \lambda p_2
-2\lambda -2)-p_2(h_{q-1}h_{r}+h_{q}h_{r-1}+2h_{q-1}h_{r-1}),
\end{eqnarray*}

\noindent where $f(q,r)=h_{r} (\lambda
h_{q-1}-h_{q-2})-h_{q-1}h_r-1$ and $h_k=\lambda
p_{k-1}(p_2-2)-p_2p_{k-2}-2p_{k-1}$. Combining with Lemma \ref{25}
and using Maple, we have
\begin{eqnarray}(x^{2}-1)^{3}x^{n+5}
P({\mathcal L}(T_4(p,q,r)),\lambda)& = & C_0(n;x)+W(p,q,r;x),\\
(x^{2}-1)^{2}x^{n+2} P({\mathcal L}(T_4(1,q,r)),\lambda)& = &
C_0^\prime (n;x)+W(1,q,r;x),
\end{eqnarray}
where $n=p+q+r+7$, $x$ satisfies $x^{2}-\lambda x+1=0$  and
\begin{eqnarray*}
C_0(n;x)& = & {x}^{
2\,n+9}-6\,{x}^{2\,n+7}-8\,{x}^{2\,n+6}+9\,{x}^{2\,n+5}+
36\,{x}^{2\,n+4}+29\,{x}^{2\,n+3}-30\,{x}^{2\,n+2}\\
& & -87\,{x}^{2\,n+1}-72\,{x}^{2\,n}
+9\,{x}^{2\,n-1}+78\,{x}^{2\,n-2
}+84\,{x}^{2\,n-3}+48\,{x}^{2\,n-4}\\
& & +15\,{x}^{2\,n-5} +2\,{x}^{2\,n
-6}-2\,{x}^{20}-15\,{x}^{19}-48\,{x}^{18}-84\,{x}^{17}-78\,{x}^{16}-9\,{x
}^{15}\\
& &
+72\,{x}^{14} +87\,{x}^{13}+30\,{x}^{12}-29\,{x}^{11}-36\,{x}^{10
}-9\,{x}^{9}+8\,{x}^{8}+6\,{x}^{7}-{x}^{5},\\
W(p,q,r;x) & = &{x}^{2\,p+7}+{x}^{2\,q+7}+{x}^{2\,r+7}+4\,{x}^{2\,p+8}
+4\,{x}^{2\,q+8
}+4\,{x}^{2\,r+8}+4\,{x}^{2\,p+9}+4\,{x}^{2\,q+9}\\
& &
+4\,{x}^{2\,r+9} -8\,{x}^{2 \,p+10}-8\,{x}^{2\,q+10}-8\,{x}^{2
\,r+10}-29\,{x}^{2\,p+11}-29\,{x}^{2\,q+11}\\
& &
-29\,{x}^{2\,r+11}
-34\,{x}^{2\,p+12}-34\,{x}^{2\,q+12}-34\,{x}^{2\,r+12}
-{x}^{2\,p+13}-{x}^{2\,q+13}\\
& &
-{x}^{2\,r+13} +52\,{x}^{2\, p+14}+52\,{x}^{2\,q+14}+52\,{x}^{2\,
r+14}+79\,{x}^{2\,p+15}
+79\,{x}^{2\,q+15}\\
& &
+79\,{x}^{2\,r+15}+58\,{x}^{2\,p+16}+58\,{x}^{2\,q+16}
+58\,{x}^{2\,r+16}+15\,{x}^{2\,p+17}+
15\,{x}^{2\,q+17}\\
& &
+15\,{x}^{2\,r+17} -12\,{x}^{2 \,p+18}
-12\,{x}^{2\,q+18}-12\,{x}^{2 \,r+18}-14\,{x}^{2\,p+19}
-14\,{x}^{2\,q+19}\\
& &
-14\,{x}^{2\,r+19} -6\,{x}^{2\,p+20}
-6\,{x}^{2\,q+20}-6\,{x}^{2\,r+20} -{x}^{2\,p+21}
 -{x}^{2\,q+21}-{x}^{2\,r+21}\\
& &
 +{x}^{2\,p+2\,
q+7}+6\,{x}^{2\,p+2\,q+8}+14\,{x}^{2\,p+2\,q+9}+12\,{x}^{2\,p+2\,q+10}
-15\,{x}^{2\,p+2\,q+11}\\
& &
-58\,{x}^{2\,p+2\,q+12}-79\,{x}^{2\,p+2\,q+13}-
52\,{x}^{2\,p+2\,q+14}+{x}^{2\,p+2\,q+15}+34\,{x}^{2\,p+2\,q+16}\\
& &
+29\,{
x}^{2\,p+2\,q+17}+8\,{x}^{2\,p+2\,q+18}-4\,{x}^{2\,p+2\,q+19}-4\,{x}^{
2\,p+2\,q+20}-{x}^{2\,p+2\,q+21}\\
& &
+{x}^{2\,p+2\,r+7}+6\,{x}^{2\,p+2\,r+8
}+14\,{x}^{2\,p+2\,r+9}+12\,{x}^{2\,p+2\,r+10}-15\,{x}^{2\,p+2\,r+11}\\
& &
-
58\,{x}^{2\,p+2\,r+12}-79\,{x}^{2\,p+2\,r+13}-52\,{x}^{2\,p+2\,r+14}+{
x}^{2\,p+2\,r+15}+34\,{x}^{2\,p+2\,r+16}\\
& &
+29\,{x}^{2\,p+2\,r+17}+8\,{x}
^{2\,p+2\,r+18}-4\,{x}^{2\,p+2\,r+19}-4\,{x}^{2\,p+2\,r+20}-{x}^{2\,p+
2\,r+21}\\
& &
+{x}^{2\,q+2\,r+7}+6\,{x}^{2\,q+2\,r+8}+14\,{x}^{2\,q+2\,r+9}+
12\,{x}^{2\,q+2\,r+10}-15\,{x}^{2\,q+2\,r+11}\\
& &
-58\,{x}^{2\,q+2\,r+12}-
79\,{x}^{2\,q+2\,r+13}-52\,{x}^{2\,q+2\,r+14}+{x}^{2\,q+2\,r+15}+34\,{
x}^{2\,q+2\,r+16}\\
& &
+29\,{x}^{2\,q+2\,r+17}+8\,{x}^{2\,q+2\,r+18}-4\,{x}^
{2\,q+2\,r+19}-4\,{x}^{2\,q+2\,r+20}-{x}^{2\,q+2\,r+21},\\
C_0^\prime (n;x)& =&{x}^{2\,n+5}-5\,{x}^{2\,n+3}-8\,{x}^{2\,n+2}
+3\,{x}^{2\, n+1}+24\,{x}^{2\,n}+28\,{x}^{2\,n-1}+2\,{x}^{2\,n-2}\\
& &
 -30\,{x
}^{2\,n-3}-36\,{x}^{2\,n-4}-20\,{x}^{2\,n-5}-10\,{x}^{2\,n-6}-15\,
{x}^{2\,n-7}-20\,{x}^{2\,n-8}\\
& &
 -15\,{x}^{2\,n-9}-6\,{x}^{2\,n-10} -{x}^{2\,n-11}
-{x}^{19}-6\,{x}^{18}-15\,{x}^{17}-20\,{x}^{16}\\
& & -15\,{x}^{15} -10\,{x}^{
14}-20\,{x}^{13}-36\,{x}^{12}-30\,{x}^{11}+2\,{x}^{10}+28\,{x}^{9}\\
& & +24 \,{x}^{8}+3\,{x}^{7}-8\,{x}^{6} -5\,{x}^{5}+{x}^{3},
\end{eqnarray*}
\begin{eqnarray*}
W(1,q,r;x) & =&-{x}^{2\,q+5} -{x}^{2\,r+5} -4\,{x}^{2\,q+6}
-4\,{x}^{2\,r+6}-6\,{x}^{2\,q+7}-6
\,{x}^{2\,r+7}-2\,{x}^{2\,q+8}\\
& &
-2\,{x}^{2\,r+8}+9\,{x}^{
2\,q+9}+9\,{x}^{2\,r+9}+20\,{x}^{2\,q+10}+20\,{x}^{2\,r+10}+25\,{x}^{2\,q+11}\\
& &
+25\,{ x}^{2\,r+11}+26\,{x}^{2\,q+12}+26\,{x}^{2\,r+12}+25\,{x}^
{2\,q+13}+25\,{x}^{2\,r+13}+20\,{x}^{2\,q+14}\\
& &
+20\,{x}^{2\,r+14}+9\,{x}^{2\,q+15}+9\,
{x}^{2\,r+15}-2\,{x}^{2\,q+16}-2\,{x}^{2\,r+16}-6\,{x}^{
2\,q+17}\\
& &
-6\,{x}^{2\,r+17}-4\,{x}^{2\,q+18}-4\,{x}^{2\,r+18}-{x}^{2\,q+19}
-{x}^{ 2\,r+19}.
\end{eqnarray*}

In view of point above, if two graphs $T_4(p,q,r)$ and
$T_4(p^\prime,q^\prime,r^\prime)$ are cospectral with respect to
Laplacian spectrum, then ${\mathcal L}(T_4(p,q,r))$ and ${\mathcal
L}(T_4(p^\prime,q^\prime,r^\prime))$ are cospectral with respect to
adjacency spectrum, hence $p+q+r=p^\prime+q^\prime+r^\prime$ and so
$W(p,q,r;x)=W(p^\prime,q^\prime,r^\prime;x)$.

\begin{lemma}\label{31}
No two non-isomorphism graphs $T_4(p,q,r)$ are cospectral with respect to Laplacian spectrum.
\end{lemma}

\noindent {\bf Proof.} Suppose that $G=T_4(p,q,r)$ and $G^\prime
=T_4(p^\prime,q^\prime,r^\prime)$ are cospectral with respect to
Laplacian spectrum. Then $G$ and $G^\prime$ have the same number of
vertices and so $p+q+r=p^\prime +q^\prime +r^\prime$. On the other
hand, by Lemma \ref{21}, ${\mathcal L}(G)$ and ${\mathcal
L}(G^\prime)$ are cospectral with respect to adjacency spectrum, so
they have the same number of closed walks of any length, especially
of length $5$. Hence ${\mathcal L}(G)$ and ${\mathcal L}(G^\prime)$
have the same number of $G_1$ in it by Lemma \ref{212}~(ii).

Clearly, for $2\leq p\leq q \leq r$, $2 \leq q^\prime \leq
r^\prime$, $2\leq r^{\prime\prime}$, $N_{{\mathcal
L}(T_4(p,q,r))}(G_1)=6$, $N_{{\mathcal
L}(T_4(1,q^\prime,r^\prime))}(G_1)=8$, $N_{{\mathcal
L}(T_4(1,1,r^{\prime\prime}))}(G_1)=10$. Hence ${\mathcal
L}(T_4(p,q,r))$, ${\mathcal L}(T_4(1,q^\prime,r^\prime))$ and
${\mathcal L}(T_4(1,1,r^{\prime\prime}))$ are non-cospectral with
each other with respect to adjacency spectrum. It follows from Lemma
\ref{21} that $T_4(p,q,r)$, $T_4(1,q^\prime,r^\prime)$ and
$T_4(1,1,r^{\prime\prime})$ are non-cospectral with each other with
respect to Laplacian spectrum.

Suppose that $G=T_4(p,q,r)$ with $p>1$. Then
$G^\prime=T_4(p^\prime,q^\prime,r^\prime)$ with $p^\prime>1$. From (3.1),
$W(p,q,r;x)=W(p^\prime,q^\prime,r^\prime;x)$. Note that $p\leq q\leq r$,
$p^\prime\leq q^\prime\leq r^\prime$ and $p+q+r=p^\prime+q^\prime+r^\prime$. It follows that $p=p^\prime$, $q=q^\prime$ and $r=r^\prime$. Therefore $G$ is isomorphic to $G^\prime$.

Let $G=T_4(1,q,r)$ with $q>1$. Then
$G^\prime=T_4(1,q^\prime,r^\prime)$ and $q^\prime>1$. By (3.2),
$W(1,q,r;x)=W(1,q^\prime,r^\prime;x)$. It follows that $q=q^\prime$
and $r=r^\prime$. Therefore $G$ is isomorphic to $G^\prime$.

If $G=T_4(1,1,r)$, then $G^\prime=T_4(1,1,r^\prime)$. It is
easy to see that $r=r^\prime$ since $G$ and $G^\prime$ have the same
number of vertices. Hence $G$ is isomorphic to $G^\prime$.

Up to now, we have completed the proof of the lemma. $\square$

\begin{lemma}\label{32} Let $G$ be a tree and  $H$ be a graph cospectral to $G$ with respect to Laplacian spectrum. If $\mu_1(G)\leq 5$, then the degree sequence of $H$
is determined by the shared spectrum.
\end{lemma}
\noindent {\bf Proof.} Let $H$ be any graph cospectral to $G$ with
respect to Laplacian spectrum. Then by Lemma \ref{211} (i) and (ii),
$H$ is also a tree. Clearly, $\Delta(G)\leq4$ by Lemmas \ref{213}.
Let $x_i$ and $y_i$ be the numbers of vertices of degree $i$ in $G$
and $H$, respectively. It follows from Lemmas \ref{211} and \ref{27}
that
\[
\left\{
\begin{array}{ll}
x_{1}+x_{2}+x_{3}+x_{4}=y_{1}+y_{2}+y_{3}+y_{4},\\
x_{1}+2x_{2}+3x_{3}+4x_{4}=y_{1}+2y_{2}+3y_{3}+4y_{4},\\
x_{1}+4x_{2}+9x_{3}+16x_{4}=y_{1}+4y_{2}+9y_{3}+16y_{4},\\
x_{1}+8x_{2}+27x_{3}+64x_{4}=y_{1}+8y_{2}+27y_{3}+64y_{4}.
\end{array}\right.
\]
It implies that $y_i=x_i$ for $i=1,2,3,4$. Hence the degree sequence
of $H$ is determined by its Laplacian spectrum. $\square$

\begin{corollary}\label{33} Let $G=T_4(p,q,r)$ and $H$ be a graph cospectral to $G$ with respect to Laplacian spectrum. Then $H$ has the same degree sequence as $G$.
\end{corollary}
\noindent {\bf Proof.} Since $G$ is a tree and $\mu_1(G)< 4.9$ by
Lemma \ref{213}, the result is  followed immediately from Lemma
\ref{32}. $\square$

\begin{lemma}\label{hlh} Let $G=T_4(p,q,r)$ and $H$ be a graph cospectral to $G$ with respect to Laplacian spectrum.
Then $H=H_1$ or $H=H_2$ (see Fig.~3) for some $l_i,k_i\geq 1$ for
$i=1,\dots,6$ and $s_j,t_j \geq 0$ for $j=1,2,3$. In particular,
${\mathcal L}(H)={\mathcal L}(H_1)$ or ${\mathcal L}(H)={\mathcal
L}(H_2)$ (see Fig.~3).
\end{lemma}
\noindent {\bf Proof.} From Lemma \ref{211} and Corollary \ref{33},
we know $H$ is a tree, having $4$ vertices of degree $3$, $6$
vertices of degree $1$ and other vertices of degree $2$. So either
all vertices of degree $3$ lie on a path or exactly $3$ vertices of
degree $3$ lie on a path and no cycle. Hence $H=H_1$ or $H=H_2$ (see
Fig.~3) for some $l_i,k_i\geq 1$ for $i=1,\dots,6$ and $s_j,t_j \geq
0$ for $j=1,2,3$. $\square$ {
\begin{figure}[htbp]
\begin{center}
\includegraphics[totalheight=9.5 cm]{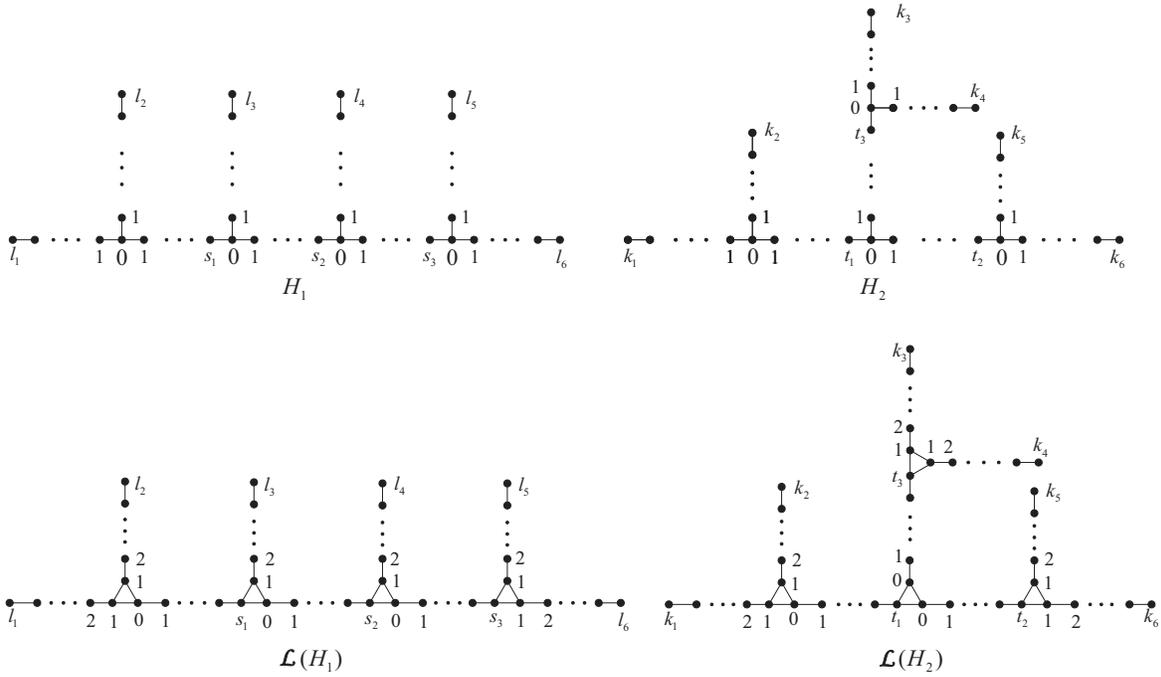}\\
 \caption{\label{HLH2} \small{The graphs $H_i$ and ${\mathcal L}(H_i)$, i=1,2, where $l_i,k_i\geq 1$ for $i=1,\dots,6$ and $s_j,t_j \geq 0$ for $j=1,2,3.$ }}
\end{center}
\end{figure}
}
\begin{lemma}\label{35} Let $G=T_4(p,q,r)$ with $p\geq 2$. Then $G$ is determined by its Laplacian spectrum.
\end{lemma}
\noindent {\bf Proof.} Let $H$ be a graph cospectral to $G$ with
respect to Laplacian spectrum. Then  ${\mathcal L}(H)$ and
${\mathcal L}(G)$ are cospectral with respect to adjacency spectrum
by Lemma \ref{21}. So ${\mathcal L}(H)$ and ${\mathcal L}(G)$ have
the same number of vertices, edges and triangles. Obviously,
$\Delta({\mathcal L}(G))=3$ and $\Delta({\mathcal L}(H))\leq 4$. Let
$y_i$ be the number of vertices of degree $i$ in ${\mathcal L}(H)$.
Note that ${\mathcal L}(G)$ has $m=p+q+r+6$ vertices, where $6$ of
them have degree $3$ and others have degree $2$. It follows from
Lemma \ref{211} that
\[
\left\{
\begin{array}{ll}
y_{1}+y_{2}+y_{3}+y_{4}=m,\\
y_{1}+2y_{2}+3y_{3}+4y_{4}=2(m+3),  \\
y_{2}+\left(3\atop2\right)y_{3}+\left(4\atop2\right)y_4=6\left(3\atop2\right)+m-6.
\end{array}\right.
\]
Solving this system of linear equation, we obtain
$(y_1,y_2,y_3,y_4)=(-y_4, m-6+3y_4,6-3y_4,y_4)$. Hence $y_1=y_4=0$
since $y_i\geq 0$ for $i=1,2,3,4$. Therefore $(y_1,y_2,y_3,y_4)=(0,
m-6,6,0)$. By Lemma \ref{hlh}, there are two cases. {
\begin{figure}[htbp]
\begin{center}
\includegraphics[totalheight=8cm]{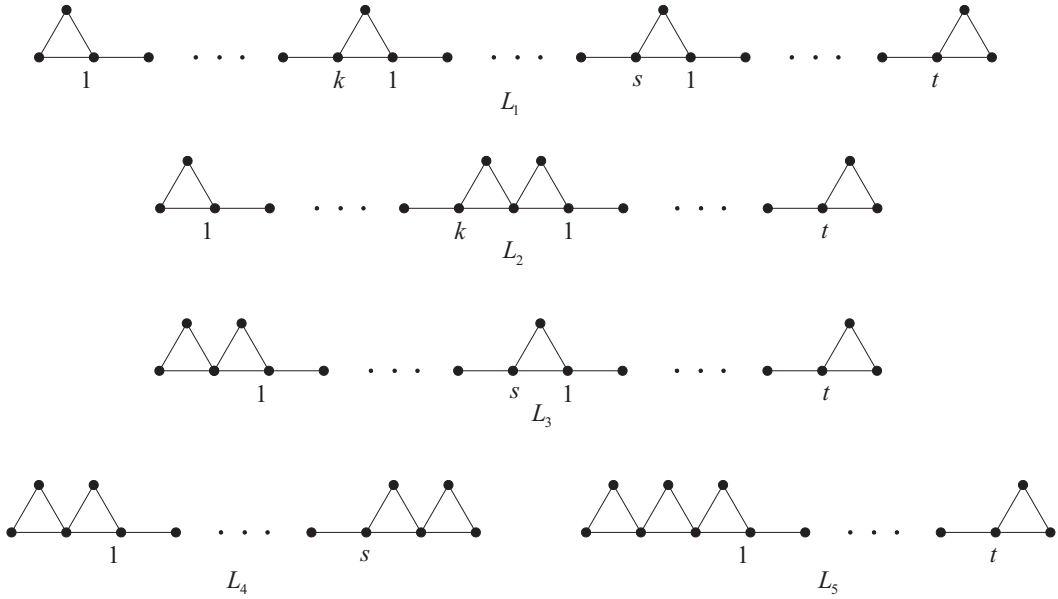}\\
 \caption{\label{l1} \small{The graphs $L_i$, $i=1,\dots,5$, where $k,s,t
 >1$.}}
\end{center}
\end{figure}
}

If ${\mathcal L}(H)={\mathcal L}(H_1)$, then $l_i=1$ and $s_j>0$ for
$i=1,\dots,6$ and $j=1,2,3$ since ${\mathcal L}(H)$ has no vertex of
degree $1$ and $4$. Hence ${\mathcal L}(H)\cong L_1$ (see Fig.~4).
Obviously, $N_{{\mathcal L}(G)}(G_1)=N_{{\mathcal L}(H)}(G_1)=6$,
$N_{{\mathcal L}(G)}(G_2)=3$, $N_{{\mathcal L}(H)}(G_2)=2$,
$N_{{\mathcal L}(G)}(G_3)=6$, $N_{{\mathcal L}(H)}(G_3)=6$ or $8$ or
$10$ or $12$, $N_{{\mathcal L}(G)}(K_3)=N_{{\mathcal L}(H)}(K_3)=4$,
$N_{{\mathcal L}(G)}(C_k)=N_{{\mathcal L}(H)}(C_k)=0$ for $k=5, 7$
and $N_{{\mathcal L}(G)}(G_i)=N_{{\mathcal L}(H)}(G_i)=0$ for
$i=4,5,6,7,8$. It follows from Lemma \ref{212} (iii) that
$N_{{\mathcal L}(G)}(7)\neq N_{{\mathcal L}(H)}(7)$. This
contradicts the fact that ${\mathcal L}(H)$ and ${\mathcal L}(G)$
are cospectral with respect to adjacency spectrum.

If ${\mathcal L}(H)={\mathcal L}(H_2)$, then $k_i=1$ and $t_j>0$ for
$i=1,\dots,6$ and $j=1,2,3$ since ${\mathcal L}(H)$ has no vertex of
degree $1$ and $4$. It implies that ${\mathcal L}(H)\cong {\mathcal
L}(T_4(p^\prime,q^\prime,r^\prime))$ for some
$p^\prime,q^\prime,r^\prime\geq 2$. Hence $H\cong
T_4(p^\prime,q^\prime,r^\prime)$ by Lemma \ref{20}. It follows from
Lemma \ref{31} that $H \cong T_4(p,q,r)=G$. $\square$

\begin{lemma}\label{36} Let $G=T_4(1,q,r)$ with $q>1$. Then $G$ is determined by its Laplacian spectrum.
\end{lemma}
\noindent {\bf Proof.} Let $H$ be a graph cospectral to $G$ with
respect to Laplacian spectrum. Then ${\mathcal L}(H)$ and ${\mathcal
L}(G)$ are cospectral with respect to adjacency spectrum by Lemma
\ref{21}. So ${\mathcal L}(H)$ and ${\mathcal L}(G)$ have the same
number of vertices, edges and triangles. Obviously,
$\Delta({\mathcal L}(G))=4$ and $\Delta({\mathcal L}(H))\leq 4$. Let
$y_i$ be the number of vertices of degree $i$ in ${\mathcal L}(H)$.
It follows from Lemma \ref{211} that
\[
\left\{
\begin{array}{ll}
y_{1}+y_{2}+y_{3}+y_{4}=m,\\
y_{1}+2y_{2}+3y_{3}+4y_{4}=2(m+3),  \\
y_{2}+\left(3\atop2\right)y_{3}+\left(4\atop2\right)y_4=\left(4\atop2\right)+4\left(3\atop2\right)+m-5,
\end{array}\right.
\]
Solving this system of linear equation, we obtain
$(y_1,y_2,y_3,y_4)=(1-y_4,m-8+3y_4,7-3y_4,y_4)$. Hence either
$y_4=0$ or $y_4=1$ since $y_i\geq 0$ for $i=1,2,3,4$.

Suppose that $y_4=0$. Then $(y_1,y_2,y_3,y_4)=(1,m-8,7,0)$, that is,
${\mathcal L}(H)$ has exactly one vertex of degree $1$, $m-8$
vertices of degree $2$, $7$ vertices of degree $3$ and no vertex of
degree $4$. Whether ${\mathcal L}(H)={\mathcal L}(H_1)$ or
${\mathcal L}(H)={\mathcal L}(H_2)$ (see Fig.~3), we always have
$N_{{\mathcal L}(H)}(G_1)= 7$, $N_{{\mathcal L}(H)}(K_3)=4$ and
$N_{{\mathcal L}(H)}(C_5)=0$. However, $N_{{\mathcal L}(G)}(G_1)=8$,
$N_{{\mathcal L}(G)}(K_3)=4$ and $N_{{\mathcal L}(G)}(C_5)=0$. It
follows from Lemma \ref{212} (ii) that $N_{{\mathcal L}(G)}(5)\neq
N_{{\mathcal L}(H)}(5)$. This contradicts the fact that ${\mathcal
L}(H)$ and ${\mathcal L}(G)$ are cospectral with respect to
adjacency spectrum.

Suppose that $y_4=1$. Then $(y_1,y_2,y_3,y_4)=(0,m-5,4,1)$. If
${\mathcal L}(H)={\mathcal L}(H_1)$, then ${\mathcal L}(H)\cong L_2$
or $L_3$ (see Fig.~4). Clearly,
\[
\begin{array}{ll}
N_{{\mathcal L}(G)}(G_1)=N_{L_2}(G_1)=N_{L_3}(G_1)=8, &
N_{{\mathcal L}(G)}(K_3)=N_{L_2}(K_3)=N_{L_3}(K_3)=4,\\
N_{{\mathcal L}(G)}(G_5)=N_{L_2}(G_5)=N_{L_3}(G_5)=2, & N_{{\mathcal L}(G)}(C_i)=N_{L_2}(C_i)=N_{L_3}(C_i)=0,\  i=5,7,\\
N_{{\mathcal L}(G)}(G_i)=N_{L_2}(G_i)=N_{L_3}(G_i)=0, & i=4,6,7,8.\\
\end{array}
\]
However,
\[
\begin{array}{lll}
N_{{\mathcal L}(G)}(G_2)=5, & N_{L_2}(G_2)=4, & N_{L_3}(G_2)=3, \\
N_{{\mathcal L}(G)}(G_3)=10\ \text{or}\ 12\ \mbox{or}\ 14, & N_{L_2}(G_3)=10\ \mbox{or}\ 12\ \mbox{or}\ 14, & N_{L_3}(G_3)=9\ \mbox{or}\ 11\ \mbox{or}\ 13.\\
\end{array}
\]
It follows from Lemma \ref{212} (iii) that $N_{{\mathcal
L}(G)}(7)\neq N_{{\mathcal L}(H)}(7)$. This contradicts the fact
that ${\mathcal L}(H)$ and ${\mathcal L}(G)$ are cospectral with
respect to adjacency spectrum.

If ${\mathcal L}(H)={\mathcal L}(H_2)$, then ${\mathcal L}(H)\cong
{\mathcal L}(T_4(1,q^\prime,r^\prime))$ for some $q^\prime,
r^\prime\geq 2$. Hence $H\cong T_4(1,q^\prime,r^\prime)$ by Lemma
\ref{20}. Therefore $H \cong T_4(1,q,r)$ by Lemma \ref{31}.
$\square$

\begin{lemma}\label{37} Let $G=T_4(1,1,r)$ with $r\geq 2$. Then $G$ is determined by its Laplacian spectrum.
\end{lemma}
\noindent {\bf Proof.} Let $H$ be a graph cospectral to $G$ with
respect to Laplacian spectrum. Then ${\mathcal L}(H)$ and ${\mathcal
L}(G)$ are cospectral with respect to adjacency spectrum by Lemma
\ref{21}. So ${\mathcal L}(H)$ and ${\mathcal L}(G)$ have the same
number of vertices, edges and triangles. Obviously,
$\Delta({\mathcal L}(G))=4$ and $\Delta({\mathcal L}(H))\leq 4$. Let
$y_i$ be the number of vertices of degree $i$ in ${\mathcal L}(H)$.
It follows from Lemma \ref{211} that
\[
\left\{
\begin{array}{ll}
y_{1}+y_{2}+y_{3}+y_{4}=m,\\
y_{1}+2y_{2}+3y_{3}+4y_{4}=2(m+3),  \\
y_{2}+\left(3\atop2\right)y_{3}+\left(4\atop2\right)y_4=2\left(4\atop2\right)+2\left(3\atop2\right)+m-4.
\end{array}\right.
\]
Solving this system of linear equation, we obtain
$(y_1,y_2,y_3,y_4)=(2-y_4, m-10+3y_4,8-3y_4,y_4)$. Hence $y_4=0$ or
$1$ or $2$ since $y_i\geq 0$ for $i=1,2,3,4$.

Suppose that $y_4=0$. Then $(y_1,y_2,y_3,y_4)=(2, m-10,8,0)$, that
is, ${\mathcal L}(H)$ has $2$ vertices of degree $1$ , $m-10$
vertices of degree $2$, $8$ vertices of degree $3$ and no vertex of
degree $4$. Whether ${\mathcal L}(H)={\mathcal L}(H_1)$ or
${\mathcal L}(H)={\mathcal L}(H_2)$, we always have $N_{{\mathcal
L}(G)}(K_3)=N_{{\mathcal L}(H)}(K_3)=4$, $N_{{\mathcal
L}(G)}(C_5)=N_{{\mathcal L}(H)}(C_5)=0$, $N_{{\mathcal L}(H)}(G_1)=
8$ and $N_{{\mathcal L}(G)}(G_1)=10$. It follows from Lemma
\ref{212} (ii) that $N_{{\mathcal L}(H)}(5)\neq N_{{\mathcal
L}(G)}(5)$. This contradicts the fact that ${\mathcal L}(H)$ and
${\mathcal L}(G)$ are cospectral with respect to adjacency spectrum.

Suppose that $y_4=1$. Then $(y_1,y_2,y_3,y_4)=(1, m-7,5,1)$, that
is, ${\mathcal L}(H)$ has $1$ vertex of degree $1$ , $m-10$ vertices
of degree $2$, $8$ vertices of degree $3$ and $1$ vertex of degree
$4$. Whether ${\mathcal L}(H)={\mathcal L}(H_1)$ or ${\mathcal
L}(H)={\mathcal L}(H_2)$, we always have $N_{{\mathcal L}(H)}(5)\neq
N_{{\mathcal L}(G)}(5)$, contradiction.

Suppose that $y_4=2$. Then $(y_1,y_2,y_3,y_4)=(0, m-4,2,2)$. If
${\mathcal L}(H)={\mathcal L}(H_1)$, then ${\mathcal L}(H)\cong L_4$
or $L_5$ (see Fig.~4). Clearly,
\[
\begin{array}{ll}
N_{{\mathcal L}(G)}(G_1)=N_{L_4}(G_1)=N_{L_5}(G_1)=10, & N_{{\mathcal L}(G)}(G_5)=N_{L_4}(G_5)= N_{L_5}(G_5)=4.\\
N_{{\mathcal L}(G)}(K_3)=N_{L_4}(K_3)=N_{ L_5}(K_3)=4, &  N_{{\mathcal L}(G)}(C_k)=N_{L_4}(C_k)=N_{L_5}(C_k)=0, k=5, 7,\\
N_{{\mathcal L}(G)}(G_i)=N_{L_4}(G_i)=N_{L_5}(G_i)=0, & i=4,6,7,8.\\
\end{array}
\]
However,
\[
\begin{array}{lll}
N_{{\mathcal L}(G)}(G_2)=8, & N_{L_4}(G_2)=4, & N_{L_5}(G_2)=6,\\
N_{{\mathcal L}(G)}(G_3)=16\ \mbox{or}\ 18, & N_{L_4}(G_3)=12\ \mbox{or}\ 14, &  N_{L_5}(G_3)=15\ \mbox{or}\ 17.\\
\end{array}
\]
It follows from Lemma \ref{212} (iii) that $N_{{\mathcal
L}(G)}(7)\neq N_{{\mathcal L}(H)}(7)$. This contradicts the fact
that ${\mathcal L}(H)$ and ${\mathcal L}(G)$ are cospectral with
respect to adjacency spectrum.

If ${\mathcal L}(H)={\mathcal L}(H_2)$, then ${\mathcal L}(H)\cong
{\mathcal L}(T_4(1,1,r^\prime))$ for some $ r^\prime\geq 2$. Hence
$H\cong T_4(1,1,r^\prime)$ by Lemma \ref{20}. Therefore $H \cong
T_4(1,1,r)$ by Lemma \ref{31}. $\square$

\begin{lemma}\label{38} Let $G=T_4(1,1,1)$. Then $G$ is determined by its Laplacian spectrum.
\end{lemma}
\noindent {\bf Proof.} Let $H$ be a graph cospectral to $G$ with
respect to Laplacian spectrum. By Lemma \ref{32}, the degree
sequence of $H$ is $(3,3,3,3,1,1,1,1,1,1)$, so $H$ is isomorphic to
a centipede graph or $T_4(1,1,1)$. By Lemma \ref{26}, the centipede
is determined by its Laplacian spectrum. Hence $H\cong T_4(1,1,1)$.
$\square$

\vskip 0.05in

Now we may give our main result in this section.
\begin{theorem}\label{39}  $T_4(p,q,r)$ is determined by its Laplacian spectrum.
\end{theorem}
\noindent {\bf Proof.} It follows from Lemmas \ref{35}, \ref{36},
\ref{37} and \ref{38}. $\square$

\vskip 0.05in

Recall from \cite{AKK} that the Laplacian eigenvalues of the
complement of a graph $G$ are completely determined by the Laplacian
eigenvalues of $G$. As a direct consequence of Theorem \ref{39}, we
have

\begin{corollary} The complement of $T_4(p,q,r)$ is determined by its Laplacian spectrum.
\end{corollary}

\section{Adjacency spectral characterization of $T_4(p,q,r)$  }

In this section, we will study the adjacency spectral
characterization of $T_4(p,q,r)$. It will be shown that there is no
two non-isomorphism graphs $T_4(p,q,r)$ are cospectral with respect
to adjacency spectrum.

Using Lemma \ref{24} with $v$ being the vertices of degree $3$, we
can compute the characteristic polynomial of $T_4(p,q,r)$ in terms
of the characteristic polynomials of paths. Put $f_r=\lambda (
p_{r+1}-p_{r-1})$ for any integer $r$. Then we have
\begin{equation*}
P(T_4(p,q,r),\lambda)=
\begin{cases}
\lambda p_2^{3}-3\lambda^{2}p_2^{2}, & if \  p= q=r=1,\\
\lambda p_2 p_2 f_r-2\lambda^{2}p_2f_r-p_2^{2}  f_{r-1}, & if \ 1= p= q< r,\\
\lambda p_2 f_q f_r-\lambda^{2}f_q f_r-p_2f_{q-1} f_r-p_2f_qf_{r-1}, & if \ 1= p< q\leq r ,\\
\lambda f_q f_p f_r-f_{q-1} f_p f_r-f_qf_{p-1} f_r-f_q f_p f_{r-1}, & if\ 2\leq p\leq q\leq r.\\
\end{cases}
\end{equation*}
Let $n=p+q+r+7$ and
$\phi(p,q,r)={x}^{n}\left({x}^{2}-1\right)^{3}P(T_4(p,q,r),\lambda)$.
By Lemma \ref{25}, we have
\begin{equation}
\phi(p,q,r))=
\begin{cases}
C_1(n;x), & if \ 1= p= q< r,\\
C_2(n;x)+U(1,q,r;x), & if \ 1= p< q\leq r ,\\
C_3(n;x)+U(p,q,r;x), & if\ 2\leq p\leq q\leq r,\\
\end{cases}
\end{equation}
where $x$ satisfies $x^{2}-\lambda x+1=0$ and
\begin{eqnarray*}
C_1(n;x)& = & 2\,{x}^{2\,n-13}-{x}^{2\,n-12}+2\,{x}^{2\,n-11}-4
\,{x}^{2\,n-9}+{x}^{2\,n-8}
-6\,{x}^{2\,n-7}+2\,{x}^{2\,n-6}\\
& & +6\,{x}^{2\,n-3}-2\,{x}^{2 \,n-2}+4\,{x}^{2\,n-1}
-{x}^{2\,n}-2\,{x}^{2\,n+1}-2\,{x}^{2\,n+3}+{x}^{2\,n+4}\\
 & &
-2\,{x}^{19}+{x}^{18}-2\,{x}^{17}+4\,{x}^{15}-{x}^{14}+6\,{x}^{13}-2\,
{x}^{12}-6\,{x}^{9}+2\,{x}^{8}-4\,{x}^{7}\\& &
+{x}^{6}+2\,{x}^{5}+2\,{x}^{3 }-{x}^{2}
\\C_2(n;x)
& = & +{x}^{2\,n-11}+{x}^{2\,n-10}-{x}^{2\, n-9}+{x}^{2\,n-8}
-2\,{x}^{2\,n-7}-3\,{x}^{2\,n-6}+{x}^{2\,n-5}\\
& & -3\,{x}^ {2\,n-4}+2\,{x}^{2\,n-3}+3\,{x}^{2\,n-2}+{x}^{2\,n-1}
+3\,{x}^{2\,n}-2\,{x}^{2\,n+1}-
{x}^{2\,n+2}\\
& &
-{x}^{2\,n+3}-{x}^{2\,n+4}+{x}^{2\,n+5}-{x}^{17}-{x}^{16}+{x}^{15}-{x}^{14}+2\,{x}^{13}+3\,{x}^{12}\\
& & -{x}^{11}+ 3\,{x}^{10}-2\,{x}^{9}-3\,{x}^{8}-{x}^{7}
-3\,{x}^{6}+2\,{x}^{5}+{x}^{4
}+{x}^{3}+{x}^{2}-x\\
C_3(n;x)& = & 2\,{x}^{2\,n-8}-{x}^{2\,n-6}
-6\,{x}^{2\,n-4}+3\,{x}^{2\,n-2}
+6\,{x}^{2\,n}-3\,{x}^{2\,n+2}-2\,{x}^{2\,n+4}\\
&
&
+{x}^{2\,n+6}-2\,{x}^{14}+{x}^{12}+6\,{x}^{10}-3\,{x}^{8}-6\,{x}^{6}+3\,{x}^{4}+2\,
{x}^{2}-1\\
U(1,q,r;x) & = &
{x}^{2\,q+4}+{x}^{2\,q+6}-3\,{x}^{2\,q+8}-3\,{x}^{2\,q+10}+3\,{x}^{2\,
q+12}+3\,{x}^{2\,q+14}-{x}^{2\,q+16}-{x}^{2\,q+18}\nonumber\\
& & +{x}^{2\,r+4}+{x}^{2
\,r+6}-3\,{x}^{2\,r+8}-3\,{x}^{2\,r+10}+3\,{x}^{2\,r+12}+3\,{x}^{2\,r+
14}-{x}^{2\,r+16}-{x}^{2\,r+18}\nonumber\\
U(p,q,r;x)&=&{x}^{2\,p+4}-3\,{x}^{2\,p+8}+3\,{x}^{2\,p+12}-{x}^{2\,p+16}\\
& & +{x}^{2\,q+
4}-3\,{x}^{2\,q+8}+3\,{x}^{2\,q+12}-{x}^{2\,q+16}\\
& &+{x}^{2\,r+4}-3\,{x}^
{2\,r+8}+3\,{x}^{2\,r+12}-{x}^{2\,r+16}+{x}^{2\,p+2\,q+4}-3\,{x}^{2\,p
+2\,q+8}\nonumber\\
& &+3\,{x}^{2\,p+2\,q+12}-{x}^{2\,p+2\,q+16}+{x}^{2\,p+2\,r+4}-3
\,{x}^{2\,p+2\,r+8}+3\,{x}^{2\,p+2\,r+12}\nonumber\\
& &-{x}^{2\,p+2\,r+16}+{x}^{2\,q
+2\,r+4}-3\,{x}^{2\,q+2\,r+8}+3\,{x}^{2\,q+2\,r+12}-{x}^{2\,q+2\,r+16}\nonumber
\end{eqnarray*}
\begin{theorem} No two non-isomorphism graphs $T_4(p,q,r)$ are cospectral with respect to adjacency spectrum.
\end{theorem}
\noindent {\bf Proof.} Suppose that $G=T_4(p,q,r)$ and
$G^\prime=T_4(p^\prime,q^\prime,r^\prime)$ are cospectral with
respect to adjacency spectrum. Then $p+q+r=p^\prime
+q^\prime+r^\prime$ and
$\phi(p,q,r)=\phi(p^\prime,q^\prime,r^\prime)$, hence
$U(p,q,r;x)=U(p^\prime,q^\prime,r^\prime;x)$. Obviously, for any
positive integers $p, q, r$ with $2\leq p\leq q\leq r$,
$\phi(1,1,r),\phi(1,q,r)$ and $\phi(p,q,r)$ are three distinct
polynomials. Therefore $T_4(1,1,r), T_4(1,q,r)$ and $T_4(p,q,r)$ are
non-cospectral with each other with respect to adjacency spectrum.

Let $G=T_4(p,q,r)$ with $2\leq p\leq q\leq r$. Then $G^\prime=T_4(p^\prime,q^\prime,r^\prime)$ with $2\leq p^\prime\leq q^\prime\leq r^\prime$ and $U(p,q,r;x)=U(p^\prime,q^\prime,r^\prime;x)$. It follows that $p=p^\prime$, $q=q^\prime$ and $r=r^\prime$. Therefore $G\cong G^\prime$.

Let $G=T_4(1,q,r)$ with $2\leq q\leq r$. Then
$G^\prime=T_4(1,q^\prime,r^\prime)$ with $2\leq q^\prime\leq
r^\prime$ and $U(1,q,r;x)=U(1,q^\prime,r^\prime;x)$. It follows that
$q=q^\prime$ and $r=r^\prime$. Therefore $G\cong G^\prime$.

Let $G=T_4(1,1,r)$ with $1\leq r$. Then $G^\prime=T_4(1,1,r^\prime)$
with $1\leq r^\prime$ and so $r=r^\prime$. Therefore $G\cong
G^\prime$.

Up to now, we have completed the proof of the theorem. $\square$

 \end{document}